\documentstyle[12pt]{article}
\def\lqq{\lq \lq }

\begin{document}

\begin{titlepage}

\begin{center}
{\LARGE Lie-algebraic discretization of differential equations
\vskip 0.5cm
\large by

\vskip 0.5cm
Yuri Smirnov$^{\star}$ and Alexander Turbiner}$^{\dagger}$
\vskip 0.5cm
Instituto de F\'isica, UNAM, Apartado Postal 20-364, 01000 Mexico D.F., Mexico
\end{center}
\vskip 0.5cm
\begin{center}
{\large ABSTRACT}
\end{center}
\vskip 0.3 cm
\begin{quote}
A certain representation for the Heisenberg algebra in finite-difference
operators
is established. The Lie-algebraic procedure of discretization of differential
equations with isospectral property is proposed. Using $sl_2$-algebra based
approach, (quasi)-exactly-solvable finite-difference equations are described.
It is shown that the operators having the Hahn, Charlier and Meixner
polynomials
 as the eigenfunctions are reproduced in present approach as some particular
cases. A discrete version of the classical orthogonal polynomials
(like Hermite, Laguerre, Legendre  and Jacobi ones) is introduced.
\end{quote}

\vfill

\noindent
$^\star$On leave of absence from the Nuclear Physics Institute,
Moscow State University, Moscow 119899, Russia\\
E-mail: smirnov@sysul2.ifisicacu.unam.mx

\noindent
$^\dagger$On leave of absence from the Institute for Theoretical
and Experimental Physics,
Moscow 117259, Russia\\
E-mail: turbiner@teorica0.ifisicacu.unam.mx  or turbiner@vxcern.cern.ch
\end{titlepage}

Discretization is one of the most powerful tools of solving continuous
theories. It appears in different forms in various physics sciences leading to
discrete versions of differential equations of classical mechanics,
lattice field theories etc.
However, once the discretization is chosen as an approach to a problem
we meet a hard problem of ambiguity -- there exist infinitely-many different
ways of discretization having the same continuous limit.
The goal of present paper is to demonstrate a certain discretization scheme for
differential equations with exceptional properties -- eigenvalues remain
unchanged (isospectrality), eigenfunctions are modified in a simple manner etc.

1.  Two operators $a$ and $b$ obeying the commutation relation
\begin{equation}
\label{e1}
[a,b] \equiv ab  -  ba \ =\ 1,
\end{equation}
with the identity operator in the r.h.s. define the Heisenberg algebra.
The Heisenberg algebra plays the central role in many branches of  theoretical
and mathematical physics. A standard representation of (1) in the action on the
real line is the coordinate-momentum one:
\begin{equation}
\label{e2}
a\ =\ {d \over dx} \ ,\ b \ =\ x \ .
\end{equation}
Our goal is to build the one-parametric representation of (1) on the line in
terms of finite-difference operators  having (2) as a limiting case.

Let us introduce the finite-difference operators ${\cal D}_{\pm}$ possessing
the property of  the translation covariance
\begin{equation}
\label{e3}
{\cal D}_+ f(x) \ =\ {f(x+\delta) - f(x) \over \delta} \equiv
{(e^{\delta {d \over dx}} -1) \over \delta} f(x)
\end{equation}
and
\begin{equation}
\label{e4}
{\cal D}_- f(x) \ =\ {f(x) - f(x - \delta) \over \delta} \equiv
{ (1 - e^{ - \delta {d \over dx}}) \over \delta} f(x)
\end{equation}
and ${\cal D}_+ \rightarrow {\cal D}_- $, once $\delta \rightarrow -\delta$.
It is worth noting that
\begin{equation}
\label{e5}
{\cal D}_+ - {\cal D}_-  \ =\ \delta {\cal D}_-{\cal D}_+
\end{equation}

Now define the so-called {\it quasi-monomial}
\begin{equation}
\label{e6}
x^{(n+1)} \ =\ x(x - \delta) (x - 2\delta) \ldots (x - n\delta) \ =\
\delta^{n+1}
{\Gamma ({x \over \delta}+1) \over
\Gamma ({x \over \delta} - n)} \ ,
\end{equation}
which in the limit $\delta \rightarrow 0$ becomes the ordinary monomial
$x^{n+1}$. The operators ${\cal D}_{\pm}$ (3)--(4) act on quasi-monomials
as follows
\begin{equation}
\label{e7}
{\cal D}_+ x^{(n)} \ =\ n x^{(n -1)} \ , \ x(1 - \delta {\cal D}_-) x^{(n)}
\ =\  x^{(n +1)} \ .
\end{equation}

So the operator ${\cal D}_+$ plays the same role in action on the
quasi-monomial
$x^{(n)} $ as the operator of differentiation in action on monomial $x^n$: it
converts $n$th quasi-monomial to $(n-1)$th with the coefficient of
proportionality equals to $n$. The operator  $x(1 - \delta {\cal D}_-) $ acts
similarly  to the operator of  multiplication $x$ on monomials: it converts
$n$th quasi-monomial to $(n+1)$th quasi-monomial. In the limit  $\delta
\rightarrow 0$ those operators
become $ {d \over dx}$  and $x $ , respectively. Straightforward calculation
shows that the commutator of the operators ${\cal D}_+$ and $x(1 - \delta
{\cal D}_-) $ equals to one ! So the operators
\[
a\ =\ {\cal D}_+ \ ,
\]
\begin{equation}
\label{e8}
b\ =\ x(1 - \delta {\cal D}_-)
\end{equation}
form one-parametric representation of the Heisenberg algebra (1). It is easy to
show that (8) belongs to the point canonical transformations.

One of possible interpretations of (8) is to identify the finite-difference
operator
$a$ with momentum in discretized form,  while  $b$ becomes the conjugated
coordinate (in a complete agreement to the limit $\delta \rightarrow 0$).

2. Consider the eigenvalue problem for some linear differential operator
\begin{equation}
\label{e9}
L [{d \over dx} , x ] \ \varphi (x) = \lambda \varphi (x)
\end{equation}
having polynomial eigenfunctions. Now let us replace in the equation (9)
the derivative ${d \over dx}$ and $x$ by the operators $a$ and $b$,
respectively, obeying the commutation relation (1). Then we arrive at
\begin{equation}
\label{e10}
L [a , b] \varphi (b) = \lambda \varphi (b)
\end{equation}
In order to make sense to (10), one should add the definition of the vacuum
$|0>$:
\begin{equation}
\label{e11}
a|0>\  = \ 0 \ .
\end{equation}
Then the equation (10) has a meaning of the operator eigenvalue problem in the
Fock space with the vacuum (11). It is evident that the eigenvalue problem
(10)--(11) has the same eigenvalues and the same eigenfunctions
(with replacement of $x \rightarrow b$) as the original problem (9).

So taking different realizations of the algebra (1) and considering a certain
element of the universal enveloping algebra of the algebra (1), we get
different
eigenvalue problems with the property of isospectrality. Below we will consider
 isospectral problems emerging from the representations (2) and (8). We name
such a way of discretization the {\it Lie-algebraic discretization}.

Now we have to specify the operators $L$ in (10) having polynomial
eigenfunctions. In \cite{t1} it was proven that $L$ has a certain amount of
polynomial eigenfunctions if and only if, $L$ is the superposition of
the element of the universal enveloping algebra of the $sl_2$-algebra taken
in the finite-dimensional irreducible representation
\[
J^+_n = b^2 a - n b
\]
\begin{equation}
\label{e12}
J^0_n = ba - {n \over 2}
\end{equation}
\[
J^-_n=a
\]
where $n$ is a non-negative integer, and the operator $B(b)a^{n+1}$,
where $B(b)$ is any operator of $b$. The dimension of this representation is
equal
to $(n+1)$ and hence $(n+1)$ eigenfunctions have a form of polynomial of degree
$n$. These operators $L$ are named {\it quasi-exactly-solvable}. Moreover,  one
can prove that $L$ possesses infinite-many polynomial eigenfunctions if and
only if $L$ is a polynomial in generators $J^0 \equiv J^0_0$, $J^- \equiv
J^-_0$
only. Those operators $L$ are named {\it exactly-solvable}.

3. In order to exploit the representation (8) let us firstly define the
vacuum $|0>$.  The condition (11) in the explicit form  is
\begin{equation}
\label{e13}
f(x+\delta) \ =\ f(x) \ .
\end{equation}
Any periodic function with the period $\delta$ is the solution of this
equation,
however, without the loss of generality we can make the choice
\begin{equation}
\label{e14}
f(x) \ =\  1 \ .
\end{equation}
With this vacuum, it is easy to see that
\begin{equation}
\label{e15}
b^n |0> = [x(1 - \delta {\cal D}_-) ]^n |0> = x^{(k)}
\end{equation}
where  the relation
\[
[x e^{ - \delta {d \over dx}}]^n \ =\ x^{(n)} e^{ -n \delta {d \over dx}} \ .
\]
was used for derivation (15).
So, starting from (10)--(11), we obtain the eigenvalue problem for the
finite-difference operator
\begin{equation}
\label{e16}
L [{\cal D}_+ , x(1 - \delta {\cal D}_-)] \tilde\varphi (x) = \lambda
\tilde\varphi (x)
\end{equation}
where the solutions $\tilde\varphi (x)$ are related to the solutions

$\varphi (x)$ of (9). To clarify this relationship, let us assume that
\begin{equation}
\label{e17}
\varphi (x) = \sum \alpha_k x^k
\end{equation}
is a certain solution of the equation (9). The transition from (9) to (10)
implies the replacement of $x$ by $b$.
Then taking into account (15) we come to the conclusion that each monomial
 in (17) should be replaced by the quasi-monomial. Finally, we get that
the corresponding solution of the equation (16) has the form
\begin{equation}
\label{e18}
\tilde\varphi (x) = \sum \alpha_k x^{(k)}
\end{equation}
The equation (16) has the same eigenvalues as the original equation (9),
while the eigenfunctions are modified by replacing the monomials by
quasi-monomials. It is evident, that such a procedure preserves the degree
of polynomial.

 Now let us proceed to concrete realization of idea of the Lie-algebraic
discretization for differential equations with polynomial solutions. It is
quite
obvious that there exist two different approaches:
\begin{itemize}
\item
In the first approach, we start from a certain differential equation with known
spectra and apply above-described Lie-algebraic discretization.
Generically, a second-order differential equation becomes a finite-difference
equation relating the unknown function in more than {\it three} points or,
in other words, more than the three-point finite-difference equations.
\item
The second approach is different: we look for the differential equations
leading to the three-point finite-difference equations under the
Lie-algebraic discretization. In other words, what are the elements of the
universal enveloping $sl_2$-algebra in the representation (8), (12) producing
the three-point finite-difference equations.
\end{itemize}
4.  {\it The first approach}. In \cite{t1} (see also \cite{t2}) it is shown
that the most general, second-order exactly-solvable differential operator
on the real line possessing infinitely-many polynomial eigenfunctions has the
form of the hypergeometrical operator
\begin{equation}
\label{e19}
E_2 ({d \over dx} , x )\ =\ - Q_2 (x) {d^2 \over dx^2} + Q_1(x) {d \over dx}
+ Q_0
\end{equation}
where $Q_k(x)$ are polynomials of the $k$th order with arbitrary coefficients:
\[
Q_2 = a_0 x^2 + a_1 x + a_2, \ Q_1= b_0 x + b_1, \ Q_0 = c_0 \ .
\]
The operator (19) is equivalent to the operator appearing from the
most general quadratic polynomial in the generators $J^0, J^-$ (see (12))
defined through the representation (2) \cite{t1}. Performing the
above-mentioned
procedure of the Lie algebraization:  $(9) \rightarrow (10) \rightarrow (16)$,
we arrive at the following exactly-solvable finite-difference operator
isospectral with (19)
\[
\tilde E_2 \equiv E_2 [{\cal D}_+ , x(1 - \delta {\cal D}_-)]\ =\
\]
\[
 -\tilde a_2 e^{2\delta {d \over dx}} + [-\tilde a_1 x +(2\tilde a_2+
\tilde b_1)]
e^{\delta {d \over dx}} +
\]
\[
 [-\tilde a_0 x(x - \delta) +(2\tilde a_1+\tilde b_0) x - (\tilde a_2+\tilde
b_1)] +
\]
\begin{equation}
\label{e20}
[2\tilde a_0 x(x - \delta) - (\tilde a_1+\tilde b_0) x] e^{-\delta {d \over
dx}} -
\tilde a_0 x (x - \delta) e^{-2\delta { d \over dx}},
\end{equation}
where $\tilde a_i ={a_i \over \delta^2}, \tilde b_i ={b_i \over \delta}$.
The corresponding eigenvalue problem has the form
\[
 -\tilde a_2 \varphi (x+2\delta) + [-\tilde a_1 x +(2\tilde a_2+\tilde b_1)]
\varphi (x+\delta) +
\]
\[
 [-\tilde a_0 x(x - \delta) +(2\tilde a_1+\tilde b_0) x - (\tilde a_2+b_1)]
\varphi (x)+
\]
\begin{equation}
\label{e21}
[2\tilde a_0 x(x - \delta) - (\tilde a_1+\tilde b_0) x] \varphi (x-\delta) -
\tilde a_0 x (x - \delta) \varphi (x-2\delta) = \lambda \varphi (x),
\end{equation}
It is worth noting that in contrary to a naive expectation, although it is
started from the second-order differential operator, the finite-difference
operator
connects to the function in five different points: $ (x +2 \delta, x + \delta,
x, x - \delta, x - 2\delta)$.

Now taking, for instance, a concrete operator (19) with the Hermite polynomials
$H_k(x)$ as the eigenfunctions $(a_0=a_1=b_1=0, a_2=-1, b_0=-2)$:
\begin{equation}
\label{e22}
h ({d \over dx}, x)= {d^2 \over dx^2} - 2x {d \over dx},
\end{equation}
or, equivalently,
\[
h \ = \ J^- J^- \ -\ 2 J^0
\]
and performing the Lie algebraization (taking into account (5)),
we arrive at the isospectral finite-difference operator
\begin{equation}
\label{e23}
\tilde h \ =\ {\cal D}_+^2 - 2 x {\cal D}_-
\end{equation}
The corresponding eigenvalue problem has the form
\[
 {1 \over \delta^2} \varphi (x+2\delta) - {2 \over \delta^2} \varphi
(x+\delta) -
 {1\over \delta} (2x  - {1\over \delta}) \varphi (x)
 + {2x \over \delta} \varphi (x-\delta) = \lambda \varphi (x),
\]
with the eigenvalues $\lambda_k=2k$ and the eigenfunctions
$H_{(k)}(x) = \sum_{i=0}^{i=k} c_i x^{(i)}$, where $k=0,1,2\ldots$ and $c_i$
are
the standard Hermite polynomials coefficients. We name these polynomials
$H_{(k)}(x)$ the {\it discrete Hermite polynomials}. In similar way one can
construct the discrete Laguerre, Legendre, Jacobi polynomials. The results
will be presented elsewhere.

Analogously, one can take the general second-order quasi-exactly-solvable
operator
\begin{equation}
\label{e24}
T_2 ({d \over dx} , x )\ =\ - P_4 (x) {d^2 \over dx^2} + P_3(x)
{d \over dx} + P_2(x)
\end{equation}
where $P_k(x)$ are polynomials of the $k$th order depending on ten free
parameters one of which is a non-negative integer, $n$  (see (12); for more
details see \cite{t1}). It leads to the following isospectral,
quasi-exactly-solvable, finite-difference operator
\newpage
\[
\tilde T_2 \equiv T_2 [{\cal D}_+ , x(1 - \delta {\cal D}_-)]\ =\
\]
\[
  \tilde \alpha_2 e^{2\delta {d \over dx}} + \tilde P_{1,1} (x) e^{\delta
{d \over dx}} +
\]
\[
 \tilde P_{0,2} (x) + x \tilde P_{-1,2} (x) e^{-\delta {d \over dx}} + x^{(2)}
\tilde P_{-2,2} (x) e^{-2\delta {d \over dx}} +
\]
\begin{equation}
\label{e25}
 x^{(3)} \tilde P_{-3,1} (x) e^{-3\delta {d \over dx}} +
\tilde \alpha_{-4} x^{(4)} e^{-4\delta {d \over dx}} ,
\end{equation}
where $\alpha$'s are parameters and $\tilde P_{j,k} (x)$ are polynomials of the
$k$th order. It is worth emphasizing that the quasi-exactly-solvable
finite-difference equation corresponding to (25) relates unknown function to
seven different points :
$ (x +2 \delta, x + \delta, x, x - \delta, x - 2\delta, x - 3\delta,
x - 4\delta)$.

5. {\it The second approach}. Standard second-order finite-difference equation
relates an unknown function at three points and has the form
\begin{equation}
\label{e26}
 A(x) \varphi (x+\delta) - B(x)\varphi (x)+ C(x) \varphi (x-\delta)
  = \lambda \varphi (x),
\end{equation}
where $A(x), B(x), C(x)$ are some functions.
One can pose a natural problem: {\it what are the most general coefficient
functions
$A(x), B(x), C(x)$ for which the equation (26) admits infinitely-many
polynomial
eigenfunctions ?} Basically, the answer is presented in \cite{t1}:
any operator with
above property can be represented as a polynomial in the generators :
\[
J^0= {x \over \delta} (1 - e^{ - \delta {d \over dx}})
\]
\begin{equation}
\label{e27}
J^-= {1 \over \delta} ( e^{  \delta {d \over dx}}-1)
\end{equation}
One can show that the most general polynomial  in the generators (27)
leading to (26) is
\begin{equation}
\label{e28}
\tilde E = A_1 J^0J^0 (J^- + {1 \over \delta}) +A_2 J^0J^- +
 A_3 J^0 + A_4 J^- + A_5
\end{equation}
and in explicit form
\newpage
\[
[{A_4 \over \delta} + {A_2 \over \delta^2} x + {A_1 \over \delta^3} x^2]
e^{  \delta {d \over dx}} +
\]
\[
+[A_5 - {A_4 \over \delta} + ({A_1 \over \delta^2} + 2 {A_2 \over \delta^2} -
{A_3 \over \delta})x - 2 {A_1 \over \delta^3} x^2]
\]
\begin{equation}
\label{e29}
+[ - ( {A_1 \over \delta^2} - {A_2 \over \delta^2} + {A_3 \over \delta}  )x+
{A_1 \over \delta^3} x^2] e^{  -\delta {d \over dx}}
\end{equation}
The spectral problem corresponding to the operator (29) is given by
\[
({A_4 \over \delta} + {A_2 \over \delta^2} x + {A_1 \over \delta^3} x^2)
f(x+\delta)+
\]
\[
-[-A_5 + {A_4 \over \delta} - ({A_1 \over \delta^2} + 2 {A_2 \over \delta^2} -
{A_3 \over \delta})x + 2 {A_1 \over \delta^3} x^2] f(x)
\]
\begin{equation}
\label{e30}
+[ - ( {A_1 \over \delta^2} - {A_2 \over \delta^2} + {A_3 \over \delta}  )x+
{A_1 \over \delta^3} x^2]f(x-\delta) = \lambda f(x) \ .
\end{equation}
In general, this spectral problem has the Hahn polynomials
$h_k^{(\alpha,\beta)} (x, N)$ as the eigenfunctions (therein we follow
the notations of \cite{nsu}). Namely,  these polynomials appear,
if $\delta = -1, A_5=0$ and
\[
A_1=-1,\ A_2= N-\beta-2,\ A_3= -\alpha-\beta-1,\ A_4=(\beta+1)(N-1) \ .
\]
Besides that, if
\[
A_1=1,\ A_2= 2-2N-\nu,\ A_3= 1-2N-\mu-\nu,\ A_4=(N+\nu-1)(N-1)
\]
the so-called analytically-continued Hanh polynomials
$\tilde h_k^{(\mu,\nu)} (x, N)$ appear, where $k=0,1,2\ldots$.

So the equation (30) corresponds to the most general exactly-solvable
finite-difference problem, while the operator (28) is the most general element
of the universal enveloping $sl_2$-algebra leading to (26). Hence the Hahn
polynomials are related to the finite-dimensional representations of a certain
cubic element of the universal enveloping $sl_2$-algebra (for a general
discussion see \cite{t2}).

Taking in (30) $\delta=1, A_5=0$ and putting
\[
A_1=0, A_2= -\mu ,\ A_3= \mu -1, A_4= \gamma \mu \ ,
\]
we reproduce the equation having the Meixner polynomials as the eigenfunctions.
Furthermore, if
\[
A_1=0, A_2= 0 ,\ A_3=  -1, A_4=  \mu \ ,
\]
the equation (30) corresponds to the equation with the Charlier polynomials as
the eigenfunctions  (for the definition of the Meixner and Charlier polynomials
 see e.g. \cite{nsu}). For a certain particular choice of the parameters,
one can reproduce the equations having Tschebyschov and Krawtchouk
polynomials as the solutions.

Among the equations (26) there also exists quasi-exactly-solvable equations
possessing a finite amount of polynomial eigenfunctions. All those equations
are classified via  the element of the universal enveloping $sl_2$-algebra
taken in the representation (8), (12)
\begin{equation}
\label{e31}
\tilde T = A_+ (J^+_n + \delta J^0_nJ^0_n) + A_1 J^0_nJ^0_n
(J^-_n + {1 \over \delta}) + A_2 J^0_nJ^-_n + A_3 J^0_n + A_4 J^-_n + A_5
\end{equation}
(cf. (28)), where $A$'s are free parameters.

In closing the authors would like to express their gratitude to the Instituto
de F\'isica, UNAM for kind hospitality extended to them. This work was
supported in part by the research CONACyT grant.

\end{document}